\documentclass[12pt, twoside, leqno]{article}
\usepackage{amsmath,amsthm}
\usepackage{amssymb,latexsym}
\usepackage{enumerate}

\usepackage{xspace}

\newtheorem{thm}{Theorem}[section]
\newtheorem{cor}[thm]{Corollary}

\newtheorem{prop}[thm]{Proposition}
\newtheorem{q}[thm]{Question}

\theoremstyle{definition}
\newtheorem{defin}[thm]{Definition}
\newtheorem{rem}{Remark}

\numberwithin{equation}{section}

\def\st{\mathchoice{:}{:}{\,:\,}{\,:\,}}

\def\size#1{\lvert#1\rvert}

\def\conc{{}^\frown\!}			
\newcommand{\forcestext}[2][{}]{\Vdash_{#1}\text{``}{#2}\text{''}}
\def\ONE{\textsc{One}\xspace}
\def\TWO{\textsc{Two}\xspace}

\def\Gdelta{\operatorname{\mathsf{G}}_{\delta}}

\begin{document}


\baselineskip=17pt


\title{Preserving the Lindel\"{o}f property under forcing extensions}
\author{Masaru Kada%
	\thanks{Supported by 
	Grant-in-Aid for Young Scientists (B) 21740080, MEXT.}}
\date{November 9, 2010}
\maketitle


\renewcommand{\thefootnote}{}

\footnote{2010 \emph{Mathematics Subject Classification}: 
	Primary 54D20; Secondary 03E40, 54G10.}
\footnote{\emph{Key words and phrases}: 
	indestructibly Lindel\"{o}f space, 
	Rothberger property, 
	P-space, 
	infinite game, 
	forcing.}

\renewcommand{\thefootnote}{\arabic{footnote}}
\setcounter{footnote}{0}


\begin{abstract}
We investigate 
preservation of the Lindel\"{o}f property 
of topological spaces 
under forcing extensions. 
We give sufficient conditions 
for a forcing notion 
to preserve several strengthenings 
of the Lindel\"{o}f property, 
such as indestructible Lindel\"{o}f property, 
the Rothberger property 
and being a Lindel\"{o}f P-space. 
\end{abstract}

\section{Introduction}\label{sec:intro}

One of several basic open problems 
about Lindel\"{o}f spaces 
asks 
what possible cardinalities are for 
Lindel\"{o}f spaces in which each point is 
$\Gdelta$. 
A number of consistency results using 
forcing or large cardinal axioms 
have been obtained. 
A fundamental issue that emerged from that work
is the following question: 
\begin{quote}
\emph{
When does a forcing extension 
preserve the Lindel\"{o}f property?} 
\end{quote}
Surprisingly, 
little seems to be 
known about this question. 
Tall \cite{Tall:lindelofpointsgd} 
introduced the notion of \emph{indestructibly Lindel\"{o}f spaces}.
A Lindel\"{o}f space is called an indestructibly Lindel\"{o}f space 
if 
it is still Lindel\"{o}f 
after 
forcing with 
any countably closed 
poset. 
Tall pointed out that 
a Lindel\"{o}f space is indestructibly Lindel\"{o}f 
if 
it 
remains 
Lindel\"{o}f 
after forcing with the 
poset 
which adjoins a Cohen subset 
of $\omega_1$ with countable conditions, 
which is just a 
particular 
instance of a countably closed 
poset. 

The class of 
spaces 
with the \emph{Rothberger property} 
is a natural and important subclass 
of the class of indestructibly 
Lindel\"{o}f spaces. 
It is of great interest 
to know which forcing notions preserve 
the Rothberger property. 
Scheepers and Tall showed 
\cite{ST:lindelof} 
that 
forcing with 
countably closed 
posets 
as well as 
the measure algebra 
preserve 
the Rothberger property.

Both 
indestructible Lindel\"{o}f property 
and 
the Rothberger property 
are nicely characterized 
in terms of 
infinite games played on topological spaces. 
Pawlikowski \cite{Pa:pointopen} 
proved 
the Rothberger property 
is equivalent to 
the non-existence of 
a winning strategy for 
the first player 
in a certain game 
played on the space. 
Scheepers and Tall 
\cite{ST:lindelof} 
proved that 
an 
indestructibly Lindel\"{o}f space 
is characterized as a space on which 
the first player 
has no 
winning strategy 
in 
a modification of 
the game which appears in Pawlikowski's theorem 
into transfinite length. 
Moreover, 
the existence of a winning strategy 
for the second player 
in the game for indestructibly Lindel\"{o}f spaces 
also 
determines 
a noteworthy class of 
spaces, 
for 
it is known that 
if a space in which each point is $\Gdelta$ 
belongs to this class 
then its cardinality is at most 
$2^{\aleph_0}$ 
\cite[Theorem~2]{ST:lindelof}.

On the other hand, 
infinite games 
played on 
posets 
have been studied by many researchers, 
mainly in connection with 
Boolean-algebraic 
or 
forcing-theoretic properties. 
One of the 
most 
significant 
results of 
those studies 
is a 
game-theoretic 
characterization of proper forcing notions. 
Also the relations between game-theoretic properties 
and various properties of forcing notions, 
such as countable closedness, 
semiproperness, 
$\alpha$-properness, 
Axiom A, 
the Sacks property 
and 
the Laver property, 
have been studied 
by 
Foreman, 
Jech, 
Veli\v{c}kovi\'c, 
Zapletal, 
Shelah, 
Ishiu, 
Kada and others. 
See 
\cite{Fo:gameba,Ishiu:alphaproper,Je:moregame,JeSh:cccba,%
	Kada:gamecd,Vel:playful,Zap:morecg} 
for further information. 

In the present paper 
we show that, 
an indestructibly Lindel\"{o}f space 
remains Lindel\"{o}f 
after forcing 
with a 
poset 
in a natural class, 
which is described using a game 
and 
larger than the class of countably closed 
posets. 
Also we show that 
the Rothberger property is preserved under 
forcing with a 
poset 
in another natural class, 
which is again described using 
a game.

We also investigate 
preservation of being a Lindel\"{o}f space 
in the class of P-spaces. 
Forcing with 
a proper 
poset 
preserves 
being a P-space. 
We will show that
a Lindel\"{o}f P-space 
remains Lindel\"{o}f 
after forcing with 
a 
poset 
in a large 
class of 
proper 
posets. 
It is an intriguing question 
if there is an example of a Lindel\"{o}f P-space which 
is 
no longer 
a Lindel\"{o}f space 
after forcing with some proper poset 
(see Question~\ref{q:destroylindelofp}).

In Section~\ref{sec:main} 
we establish a general preservation theorem 
stated in terms of games, 
which yields all the preservation results mentioned above. 
We will prove this 
theorem 
by 
pursuing 
moves in two games played in parallel, 
one is played on a topological space 
and the other on a poset.

The general investigation 
of preservation of the Lindel\"{o}f property 
and its strengthenings 
under forcing extension 
has internal appeal, 
but the results may be useful 
in obtaining consistency results about a number of 
other basic open problems about Lindel\"{o}f spaces.

\section{Preliminaries}\label{sec:pre}

For a 
poset 
$\mathbb{P}$, 
an ordinal $\alpha$ 
and a cardinal $\kappa$, 
the 
\emph{cut-and-choose game 
$\operatorname{\mathsf{CG}}^{<\alpha}({<}\kappa)$ 
on $\mathbb{P}$} is defined as follows. 
The game is played by two players \ONE and \TWO 
for $\alpha$ innings. 
In the beginning 
\ONE chooses $p\in\mathbb{P}$. 
In each inning $\beta<\alpha$, 
\ONE chooses a $\mathbb{P}$-name $\dot{\eta}_{\beta}$ for an ordinal, 
and 
then 
\TWO chooses a set $C_{\beta}$ of ordinals 
with $\size{C_{\beta}}<\kappa$. 
\TWO wins in this game 
if for every $\gamma<\alpha$ there is 
$q_{\gamma}\in\mathbb{P}$ 
such that 
$q_{\gamma}\leq p$ 
and 
$q_{\gamma}\forcestext[\mathbb{P}]{\forall \beta<\gamma\,
	(\dot{\eta}_{\beta}\in \check{C}_{\beta})}$. 
Note that, for \TWO to win, 
it is not 
required 
to find 
$q\in\mathbb{P}$ 
such that 
$q\leq p$ 
and 
$q\forcestext[\mathbb{P}]{\forall \beta<\alpha\,
	(\dot{\eta}_{\beta}\in \check{C}_{\beta})}$. 
Sometimes we preliminarily fix 
\ONE's beginning move $p\in\mathbb{P}$ 
and then 
start 
the innings; 
in such a case we call it 
the \emph{game $\operatorname{\mathsf{CG}}^{<\alpha}({<}\kappa)$ 
on $\mathbb{P}$ below $p$}. 
If $\alpha=\gamma+1$, 
we 
write $\operatorname{\mathsf{CG}}^{\gamma}({<}\kappa)$ 
instead of $\operatorname{\mathsf{CG}}^{<\gamma+1}({<}\kappa)$. 
Also, we 
write 
$\operatorname{\mathsf{CG}}^{<\alpha}(\lambda)$ 
instead of 
$\operatorname{\mathsf{CG}}^{<\alpha}({<}\lambda^{+})$.

The following theorem is well-known 
\cite{Je:moregame}. 

\begin{thm}\label{thm:countablechoiceproper}
For a 
forcing notion 
$\mathbb{P}$, 
if \TWO has a winning strategy 
in $\operatorname{\mathsf{CG}}^{\omega}(\aleph_0)$ 
on $\mathbb{P}$, 
then $\mathbb{P}$ is proper. 
\end{thm}


We say 
a forcing notion $\mathbb{P}$ is \emph{$\omega^\omega$-bounding} 
if, 
for $p\in\mathbb{P}$ 
and a $\mathbb{P}$-name $\dot{f}$ for an element of $\omega^\omega$, 
there are $q\in\mathbb{P}$ and $g\in\omega^\omega$ 
such that $q\leq p$ and 
$q\forcestext[\mathbb{P}]{%
	\forall n<\omega\,(\dot{f}(n)\leq g(n))}$. 
The following fact is easily checked.

\begin{thm}\label{thm:finitechoicebounding}
For a 
forcing notion 
$\mathbb{P}$, 
if \TWO has a winning strategy 
in $\operatorname{\mathsf{CG}}^{\omega}({<}\aleph_0)$ 
on $\mathbb{P}$, 
then $\mathbb{P}$ is $\omega^\omega$-bounding. 
\end{thm}

\begin{rem}
Although the converse of 
Theorem~\ref{thm:finitechoicebounding}
does not hold, 
most well-known proper $\omega^\omega$-bounding 
forcing notions, 
such as Sacks forcing, Silver forcing and the measure algebra, 
are the ones 
on which 
\TWO has a winning strategy 
in $\operatorname{\mathsf{CG}}^{\omega}({<}\aleph_0)$. 
\end{rem}

A 
poset 
$\mathbb{P}$ 
is \emph{${<}\alpha$-closed} 
if any descending sequence in $\mathbb{P}$ 
of length 
less 
than $\alpha$ 
has a lower bound in $\mathbb{P}$. 
If 
$\mathbb{P}$ is ${<}\alpha$-closed, 
then 
obviously 
\TWO has a winning strategy 
in 
$\operatorname{\mathsf{CG}}^{<\alpha}(1)$ 
on $\mathbb{P}$.

\begin{rem}
The game 
	$\operatorname{\mathsf{CG}}^{<\alpha}(1)$ 
	on a poset $\mathbb{P}$ 
is closely related to 
the 
\emph{strategic closure} of $\mathbb{P}$. 
A poset $\mathbb{P}$ 
is 
\emph{${<}\alpha$-strategically closed} 
if 
the second player has a winning strategy 
in 
the 
\emph{descending chain game} 
on $\mathbb{P}$ of length $\alpha$, 
which is a generalization of a usual Banach--Mazur game 
into transfinite length 
but 
the second player has the initiative 
in each limit inning 
(see \cite{Fo:gameba} or \cite{IY:games} for a precise definition). 
For an ordinal $\alpha$ 
which is either a limit 
or the successor of a limit, 
$\mathbb{P}$ is ${<}\alpha$-strategically closed 
if and only if 
\TWO has a winning strategy 
in $\operatorname{\mathsf{CG}}^{<\alpha}(1)$ 
on $\mathbb{P}$ 
(it was proved in the case $\alpha=\omega+1$ 
by 
Jech and 
Veli\v{c}kovi\'{c} 
\cite{Je:multi,Vel:playful}, 
and in a general case by Ishiu 
in his unpublished paper~\cite{Ishiu:transgame}). 

It is unprovable in ZFC 
that 
if \TWO has a winning strategy 
in $\operatorname{\mathsf{CG}}^{<\omega_1}(1)$ 
on $\mathbb{P}$ 
then 
$\mathbb{P}$ is ${<}\omega_1$-closed, 
for the following reason: 
It is known that 
$\mathbb{P}$ is ${<}\omega_1$-strategically closed 
if and only if 
$\mathbb{P}$ is ${<}(\omega+1)$-strategically closed 
(see \cite{Vel:square} or \cite{IY:games}), 
and 
it is 
known to be 
unprovable in ZFC 
that 
if 
$\mathbb{P}$ is ${<}(\omega+1)$-strategically closed 
then 
$\mathbb{P}$ is ${<}\omega_1$-closed 
(due to Jech and Shelah \cite{JeSh:cccba}). 
\end{rem}

Here we list properties of a forcing notion $\mathbb{P}$ 
which are relevant to the results in this paper. 
$\operatorname{Fn}(\omega_1,2,\omega_1)$ 
denotes 
the 
poset 
which adjoins a Cohen subset of $\omega_1$ 
with countable conditions 
\cite{Ku:set}. 
The list is ordered \emph{stronger to weaker}. 

\begin{enumerate}
\item 
	$\mathbb{P}=\operatorname{Fn}(\omega_1,2,\omega_1)$. 
\item 
	$\mathbb{P}$ is ${<}\omega_1$-closed. 
\item 
	\TWO has a winning strategy 
	in $\operatorname{\mathsf{CG}}^{<\omega_1}(1)$ 
	on $\mathbb{P}$. 
\item 
	\TWO has a winning strategy 
	in $\operatorname{\mathsf{CG}}^{<\omega_1}({<}\aleph_0)$ 
	on $\mathbb{P}$. 
\item 
	\TWO has a winning strategy 
	in $\operatorname{\mathsf{CG}}^{\omega}({<}\aleph_0)$ 
	on $\mathbb{P}$. 
\item 
	\TWO has a winning strategy 
	in $\operatorname{\mathsf{CG}}^{\omega}(\aleph_0)$ 
	on $\mathbb{P}$. 
\item $\mathbb{P}$ is proper. 
\end{enumerate}

Now we turn to 
the games played on topological spaces. 

For a topological space $(X,\tau)$ 
and an ordinal $\alpha$, 
the 
\emph{game 
$\operatorname{\mathsf{G}}_1^{<\alpha}(\mathcal{O},\mathcal{O})$ 
on $(X,\tau)$} 
is played by two players \ONE and \TWO 
for $\alpha$ innings as follows. 
In the inning $\beta<\alpha$, 
\ONE chooses an open cover $\mathcal{U}_\beta$ of $X$ 
and then 
\TWO chooses 
$H_\beta\in\mathcal{U}_\beta$. 
\TWO wins in this game 
if there is $\gamma<\alpha$ 
such that, 
$\{H_\beta\st \beta<\gamma\}$ covers $X$. 
Note that \TWO does not win if 
just $\{H_\beta\st \beta<\alpha\}$ covers $X$.
If $\alpha=\gamma+1$, 
We 
write $\operatorname{\mathsf{G}}_1^{\gamma}(\mathcal{O},\mathcal{O})$ 
instead of $\operatorname{\mathsf{G}}_1^{<\gamma+1}(\mathcal{O},\mathcal{O})$.

We say 
a space $(X,\tau)$ 
\emph{has the Rothberger property}
if, 
for every sequence $\langle\mathcal{U}_n\st n<\omega\rangle$ 
of open covers of $X$ 
there is 
an open cover  $\{U_n\st n<\omega\}$ of $X$ 
such that 
$U_n\in\mathcal{U}_n$ for all $n<\omega$. 
It is easy to see that, 
if \ONE does not have a winning strategy 
in the game 
$\operatorname{\mathsf{G}}_1^{\omega}(\mathcal{O},\mathcal{O})$ 
on $(X,\tau)$, 
then $(X,\tau)$ has the Rothberger property. 
The following theorem, 
which is due to Pawlikowski \cite{Pa:pointopen}, 
tells us that the converse also holds. 

\begin{thm}\label{thm:rothbergergame}
A space $(X,\tau)$ has the Rothberger property 
if and only if 
\ONE 
does not have a winning strategy 
in the game 
$\operatorname{\mathsf{G}}_1^{\omega}(\mathcal{O},\mathcal{O})$ 
on $(X,\tau)$. 
\end{thm}

For a 
space $(X,\tau)$ 
and a forcing notion $\mathbb{P}$, 
we let $\tau^{\mathbb{P}}$ denote 
a $\mathbb{P}$-name representing the topology on $X$ 
generated by $\tau$ in a generic extension by $\mathbb{P}$.

We say 
a forcing notion $\mathbb{P}$ 
\emph{destroys} 
a Lindel\"{o}f space $(X,\tau)$ 
if we have 
\[
\forcestext[\mathbb{P}]{
	(\check{X},\tau^{\mathbb{P}})\text{ is not Lindel\"{o}f}}.
\]
A Lindel\"{o}f space $(X,\tau)$ is called an 
\emph{indestructibly Lindel\"{o}f} space 
if 
$(X,\tau)$ 
is not 
destroyed
by 
any ${<}\omega_1$-closed 
poset. 

The equivalence $(1)\Leftrightarrow(2)$ in the following theorem is 
due to Scheepers and Tall \cite[Theorem~1]{ST:lindelof}. 
The equivalence $(2)\Leftrightarrow(3)$ is easily checked. 

\begin{thm}\label{thm:indlindelofgame}
For a space $(X,\tau)$ the following are equivalent. 
\begin{enumerate}
\item[\textup{(1)}] 
	$(X,\tau)$ is an indestructibly Lindel\"{o}f space. 
\item[\textup{(2)}] 
	$(X,\tau)$ 
	is a Lindel\"{o}f space and 
	\ONE does not have a winning strategy in 
	$\operatorname{\mathsf{G}}_1^{\omega_1}(\mathcal{O},\mathcal{O})$ 
	on $(X,\tau)$. 
\item[\textup{(3)}] 
	\ONE does not have a winning strategy in 
	$\operatorname{\mathsf{G}}_1^{<\omega_1}(\mathcal{O},\mathcal{O})$ 
	on $(X,\tau)$. 
\end{enumerate}
\end{thm}

As a consequence of Theorem~\ref{thm:rothbergergame} 
and Theorem~\ref{thm:indlindelofgame}, 
we see 
the following fact \cite[Corollary~10]{ST:lindelof}. 

\begin{cor}
A space with the Rothberger property 
is an indestructibly Lindel\"{o}f space. 
\end{cor}

We say 
$X$ 
is a \emph{P-space} 
if 
every $\Gdelta$-set in 
$X$ 
is an open set. 
It is known that 
a Lindel\"{o}f P-space 
has the Rothberger property 
(due to Galvin; 
see the following remark). 

\begin{rem}
An open cover $\mathcal{U}$ of a space $X$ 
is 
an \emph{$\omega$-cover} 
if 
$X\notin\mathcal{U}$ and 
for every finite set $F\subseteq X$ 
there is a $U\in\mathcal{U}$ 
with $F\subseteq U$. 
An open cover $\mathcal{U}$ of 
$X$ 
is 
a \emph{$\gamma$-cover} 
if $\mathcal{U}$ is infinite and 
any infinite subset of $\mathcal{U}$ covers $X$. 
A space $X$ is called a \emph{$\gamma$-space} 
if, 
for 
every 
sequence 
$\langle\mathcal{U}_n\st n<\omega\rangle$ 
of $\omega$-covers of $X$ 
there is 
a $\gamma$-cover 
$\{U_n\st n<\omega\}$ 
of $X$ 
such that 
$U_n\in\mathcal{U}_n$ for all $n<\omega$. 
It is known that a $\gamma$-space has the Rothberger property. 
Galvin 
proved 
that 
a Lindel\"{o}f P-space 
is a $\gamma$-space 
(see \cite{GN:cx}, 
\cite[Theorem~47]{ST:lindelof}). 
\end{rem}

Here we list properties of a topological space $X=(X,\tau)$ 
which are relevant to the results in this paper. 
The list is ordered \emph{weaker to stronger}. 

\begin{enumerate}
\item $X$ is a Lindel\"{o}f space. 
\item $X$ is an indestructibly Lindel\"{o}f space 
	(equivalently, 
	\ONE does not have a winning strategy in 
	$\operatorname{\mathsf{G}}_1^{<\omega_1}(\mathcal{O},\mathcal{O})$
	on $X$).
\item $X$ has the Rothberger property 
	(equivalently, 
	\ONE does not have a winning strategy in 
	$\operatorname{\mathsf{G}}_1^{\omega}(\mathcal{O},\mathcal{O})$
	on $X$).
\item $X$ is a Lindel\"{o}f P-space. 
\end{enumerate}

\begin{rem}
Here 
we state facts about the topological property 
``\TWO has a winning strategy in 
	$\operatorname{\mathsf{G}}_1^{<\omega_1}(\mathcal{O},\mathcal{O})$
	on $X$'', 
which does not fit in the above list. 
Clearly, 
if \TWO has a winning strategy in 
	$\operatorname{\mathsf{G}}_1^{<\omega_1}(\mathcal{O},\mathcal{O})$
	on $X$, 
then $X$ is an indestructibly Lindel\"{o}f space. 
Daniels and Gruenhage 
\cite{DG:pointopentype} 
showed 
that 
	if $X$ is a hereditarily Lindel\"{o}f space 
	then 
	\TWO has a winning strategy in 
	$\operatorname{\mathsf{G}}_1^{<\omega_1}(\mathcal{O},\mathcal{O})$
	on 
	$X$. 
The real line $\mathbb{R}$ is a hereditarily Lindel\"{o}f space 
and so 
\TWO has a winning strategy in 
	$\operatorname{\mathsf{G}}_1^{<\omega_1}(\mathcal{O},\mathcal{O})$
	on $\mathbb{R}$, 
whereas 
$\mathbb{R}$ does not have the Rothberger property. 
On the other hand, 
using results due to Scheepers and Tall 
\cite[Theorem 2 and Example 3]{ST:lindelof} 
we can see that, 
it is consistent 
with ZFC 
that 
there is a space 
with the Rothberger property 
on which 
\TWO does not have a winning strategy 
in 
$\operatorname{\mathsf{G}}_1^{<\omega_1}(\mathcal{O},\mathcal{O})$.
\end{rem}

\section{Proper forcing preserves P-spaces}

We prove that a P-space remains a P-space 
after forcing with a proper 
forcing notion. 

\begin{prop}\label{prop:preservepspace}
Suppose that a space $(X,\tau)$ is a P-space 
and $\mathbb{P}$ is a proper forcing notion. 
Then 
$\forcestext[\mathbb{P}]{
	(\check{X},\tau^{\mathbb{P}})\text{ is a P-space}}$. 
\end{prop}

\begin{proof}
Fix a countable set $\{\dot{G}_n\st n<\omega\}$ 
of $\mathbb{P}$-names such that 
$
	\forcestext[\mathbb{P}]{
	\dot{G}_n\in\tau^{\mathbb{P}}
	}
$
for all $n<\omega$. 
For each $n$, 
take a $\mathbb{P}$-name $\dot{\mathcal{T}}_n$ 
such that 
\[
\forcestext[\mathbb{P}]{
	\dot{\mathcal{T}}_n\subseteq\check{\tau}\text{ and }
	\dot{G}_n=\bigcup\dot{\mathcal{T}}_n}. 
\]

We are going to prove the following sentence. 
\[
	\forcestext[\mathbb{P}]{
		\forall x\in \check{X}\,
		\left[
		x\in\bigcap_{n<\omega}\dot{G}_n
		\to
		\exists H\in\check{\tau}\,(x\in H
		\text{ and }H\subseteq\bigcap_{n<\omega}\dot{G}_n)
		\right]
		}, 
\]
which implies that 
$\forcestext[\mathbb{P}]{
	\bigcap_{n<\omega}\dot{G}_n\in\tau^{\mathbb{P}}}
$. 
It suffices to show that, 
for 
$x\in X$ and $p\in \mathbb{P}$, 
if 
$p\forcestext[\mathbb{P}]{
	\check{x}\in\bigcap_{n<\omega}\dot{G}_n}$,
then there are $q\leq p$ and $H\in\tau$ such that 
$x\in H$ and 
$q\forcestext[\mathbb{P}]{
	\check{H}\subseteq\bigcap_{n<\omega}\dot{G}_n}$.

Fix $x\in X$, $p\in\mathbb{P}$ 
and assume 
$p\forcestext{\check{x}\in\bigcap_{n<\omega}\dot{G}_n}$. 
For each $n<\omega$, 
since we have 
$
	p\forcestext[\mathbb{P}]{
		\check{x}\in\dot{G}_n\text{ and }
		\dot{G}_n=\bigcup\dot{\mathcal{T}}_n} 
$, 
we can take a $\mathbb{P}$-name $\dot{T}_n$ 
such that 
$
	p\forcestext[\mathbb{P}]{
		\check{x}\in\dot{T}_n\text{ and }
		\dot{T}_n\in\dot{\mathcal{T}}_n}
$. 
Note that we have 
$p\forcestext[\mathbb{P}]{
	\dot{T}_n\in\check{\tau}
	\text{ and }\dot{T}_n\subseteq\dot{G}_n}$. 

By the properness of $\mathbb{P}$, 
we can choose $q\leq p$ 
and a countable set $\mathcal{C}\subseteq\tau$ 
so that 
$
	q\forcestext[\mathbb{P}]{
	\{\dot{T}_n\st n<\omega\}\subseteq\check{\mathcal{C}}}
$. 
Note that 
$q\forcestext[\mathbb{P}]{
	\forall n<\omega\,(\check{x}\in\dot{T}_n)}$. 
Let $H=\bigcap\{T\in\mathcal{C}\st x\in T\}$. 
Then $x\in H$ and, 
since $(X,\tau)$ is a P-space, 
$H\in\tau$ holds. 
Now we have 
\[
	q\forcestext[\mathbb{P}]{
		\check{H}=\bigcap{\{T\in\check{\mathcal{C}}\st \check{x}\in T\}}
		\subseteq\bigcap_{n<\omega}\dot{T}_n
		\subseteq\bigcap_{n<\omega}\dot{G}_n},
\]
which concludes the proof. 
\end{proof}

\section{The main result}\label{sec:main}

In this section, 
we give a sufficient condition 
for a topological space $(X,\tau)$ 
and a forcing notion $\mathbb{P}$ 
to 
keep 
$(\check{X},\tau^{\mathbb{P}})$ 
having 
a certain 
game-theoretic property 
in the forcing extension.

\begin{defin}
For a topological space $X=(X,\tau)$, 
define a cardinal $p(X)$ 
by letting 
$
p(X)
=\aleph_0+
\min\left(
	\{\size{\mathcal{G}}\st
	\mathcal{G}\subseteq\tau
	\text{ and }
	\bigcap\mathcal{G}\notin\tau
	\}
	\cup\{\size{\tau}^{+}\}\right)
$. 
\end{defin}

Note that $X$ is a P-space if and only if $p(X)\geq\aleph_1$. 

\begin{thm}\label{thm:main}
Let $(X,\tau)$ be a 
topological 
space, 
$\mathbb{P}$ a forcing notion, 
$\alpha$ an ordinal 
and $\kappa=p(X)$. 
If 
\begin{enumerate}
\item \ONE 
	does not have a winning strategy in 
	$\operatorname{\mathsf{G}}_1^{<\alpha}(\mathcal{O},\mathcal{O})$ 
	on $(X,\tau)$, and  
\item 
	\TWO has a winning strategy in 
	$\operatorname{\mathsf{CG}}^{<\alpha}({<}\kappa)$ on $\mathbb{P}$, 
\end{enumerate}
then 
\[
	\forcestext[\mathbb{P}]{%
		\text{\ONE 
		does not have a 
		winning strategy in }
		\operatorname{\mathsf{G}}_1^{<\alpha}(\mathcal{O},\mathcal{O})
		\text{ on }
		(\check{X},\tau^{\mathbb{P}})
	}.
\]
\end{thm}

\begin{proof}
Fix an enumeration of $\tau$, 
say $\tau=\{T_\xi\st\xi<\theta\}$ 
for some cardinal $\theta$. 

Suppose that 
$\dot{\sigma}$ 
is a $\mathbb{P}$-name such that 
\[
	\forcestext[\mathbb{P}]{
		\dot{\sigma}
		\text{ is a strategy for \ONE in }
		\operatorname{\mathsf{G}}_1^{<\alpha}(\mathcal{O},\mathcal{O})
		\text{ on }
		(\check{X},\tau^{\mathbb{P}})
	}. 
\]
Without loss of generality we may assume that, 
it is forced that 
the strategy $\dot{\sigma}$ 
suggests 
only open covers 
which consist of elements of $\tau$, 
since 
$\tau$ is a base of $\tau^{\mathbb{P}}$ 
in a generic extension, 
and taking refinements will not help \TWO win easier. 
Under this assumption, 
a sequence of initial moves for \TWO, 
played in a generic extension,  
against the strategy $\dot{\sigma}$ 
will be 
described 
in a form 
$\langle\check{T}_{\dot{\xi}_\beta}\st\beta<\delta\rangle$, 
where $\delta<\alpha$ and 
each $\dot{\xi}_\beta$ is a $\mathbb{P}$-name for an ordinal.

We will prove the following statement: 
For any $p\in\mathbb{P}$, 
there are $q\leq p$, 
a sequence $\langle\dot{\xi}_\beta\st\beta<\alpha\rangle$ 
of $\mathbb{P}$-names for ordinals 
and 
$\gamma<\alpha$ 
such that 
\[
	q\forcestext[\mathbb{P}]{
	\forall\delta<
	\gamma
	\,
	\big(\check{T}_{\dot{\xi}_{\delta}}\in 
	\dot{\sigma}(\langle\check{T}_{\dot{\xi}_\beta}\st\beta<\delta\rangle)
	\big)
	\text{ and }
	\bigcup\{\check{T}_{\dot{\xi}_{\delta}}\st \delta<\gamma\}
	=\check{X}
	}. 
\]
This means that 
$\langle\check{T}_{\dot{\xi}_\beta}\st\beta<\alpha\rangle$ 
describes 
winning moves for \TWO 
against the given 
strategy 
$\dot{\sigma}$ 
for \ONE 
in a generic extension.

Fix $p\in\mathbb{P}$. 
By the assumption, 
\TWO has a winning strategy $\rho$ in 
the game 
	$\operatorname{\mathsf{CG}}^{<\alpha}({<}\kappa)$ on $\mathbb{P}$ 
	below $p$. 

We are going to define a strategy $\Sigma$ for \ONE in the game 
$\operatorname{\mathsf{G}}_1^{<\alpha}(\mathcal{O},\mathcal{O})$ 
on $(X,\tau)$, 
which cannot be a winning strategy by the assumption. 

We construct $\Sigma$ by induction on $\delta<\alpha$. 
As an additional induction hypothesis 
we assume that, with each sequence 
$\langle H_\beta\st\beta<\delta\rangle$ 
describing \TWO's 
possible 
initial moves 
against $\Sigma$ 
before the inning $\delta$, 
a sequence $\langle\dot{\xi}_\beta\st\beta<\delta\rangle$ 
of $\mathbb{P}$-names for ordinals is associated.  
We will define \ONE's move 
$\Sigma\left(\langle H_\beta\st\beta<\delta\rangle\right)$ 
in the inning $\delta$, 
and associate a $\mathbb{P}$-name $\dot{\xi}_\delta$ with 
\TWO's response $H_\delta$.

Let $\dot{\mathcal{U}}^\delta$ be a $\mathbb{P}$-name 
such that 
$\forcestext[\mathbb{P}]{
	\dot{\mathcal{U}}^\delta
	=\dot{\sigma}(\langle \check{T}_{\dot{\xi}_\beta}
	\st\beta<\delta\rangle)}$. 
Since 
we have 
$\forcestext[\mathbb{P}]{
	\dot{\mathcal{U}}^\delta\subseteq\check{\tau}
	\text{ and }
	\bigcup\dot{\mathcal{U}}^\delta=\check{X}
	}$, 
for each $x\in X$ 
we can take 
a $\mathbb{P}$-name $\dot{\eta}_x^\delta$ for an ordinal 
so that 
\[
	\forcestext[\mathbb{P}]{
	\check{x}\in \check{T}_{\dot{\eta}_x^\delta}
	\text{ and }
	\check{T}_{\dot{\eta}_x^\delta}\in\dot{\mathcal{U}}^\delta
	}. 
\]
For each $x\in X$, 
let 
$F_x
	=F\left(\langle\dot{\xi}_\beta\st\beta<\delta\rangle,x\right)
	=\rho\left(\langle\dot{\xi}_\beta\st\beta<\delta\rangle
	\conc\langle\dot{\eta}_x^\delta\rangle\right)$, 
and 
$G_x
	=G\left(\langle\dot{\xi}_\beta\st\beta<\delta\rangle,x\right)
	=\bigcap\{T_\xi\st\xi\in F_x\text{ and }x\in T_\xi\}$. 
Note that 
since $\size{F_x}<\kappa=p(X)$ 
and by the definition of $p(X)$, 
$G_x$ is an open set containing $x$. 
Now let 
\[
	\Sigma\left(\langle H_\beta\st\beta<\delta\rangle\right)
	=\{G_x\st x\in X\}. 
\]
Suppose that \TWO picks $H_\delta$ from the cover $\{G_x\st x\in X\}$ 
as a move in the inning $\delta$. 
We pick $x_\delta\in X$ such that $H_\delta=G_{x_\delta}$, 
and let $\dot{\xi}_\delta=\dot{\eta}_{x_\delta}^\delta$. 
This completes the induction step at $\delta$. 

Since $\Sigma$ is not a winning strategy, 
we can find a sequence 
$\langle H_\beta\st\beta<\alpha\rangle$
which describes \TWO's winning moves 
against $\Sigma$, 
and the associated sequence  
$\langle\dot{\xi}_\beta\st\beta<\alpha\rangle$ 
of $\mathbb{P}$-names of ordinals. 
For each $\delta<\alpha$, 
let $F_\delta=\rho(\langle\dot{\xi}_\beta\st\beta\leq\delta\rangle)$. 
Find $\gamma<\alpha$ 
such that 
$\{H_\beta\st\beta<\gamma\}$ covers $X$. 
Since $\rho$ is a winning strategy for \TWO 
in the game 
$\operatorname{\mathsf{CG}}^{<\alpha}({<}\kappa)$ on $\mathbb{P}$ 
below $p$,
we can find 
$q\leq p$ 
such that 
$	q\forcestext[\mathbb{P}]{
	\forall\delta<\gamma\,
	(\dot{\xi}_\delta\in \check{F}_\delta)
	}$. 

Fix $\delta<\gamma$. 
By the construction of 
the sequence 
$\langle\dot{\xi}_\beta\st\beta<\alpha\rangle$, 
we have 
$H_\delta
	=G\left(\langle\dot{\xi}_\beta\st\beta<\delta\rangle,x_\delta\right)
$
and 
$\dot{\xi}_\delta=\dot{\eta}_{x_\delta}^\delta$ 
for a suitable $x_\delta\in X$. 
Note that 
$F_\delta
	=\rho(\langle\dot{\xi}_\beta\st\beta\leq\delta\rangle)
	=F\left(\langle\dot{\xi}_\beta\st\beta<\delta\rangle,x_\delta\right)$ 
and so 
$H_\delta
	=\bigcap\{T_\xi\st\xi\in F_\delta\text{ and }x_\delta\in T_\xi\}
$. 
Also 
\[
	\forcestext[\mathbb{P}]{
	\check{x}_\delta\in\check{T}_{\dot{\xi}_\delta}
	\text{ and }
	\check{T}_{\dot{\xi}_\delta}\in
	\dot{\sigma}
	(\langle\check{T}_{\dot{\xi}_\beta}\st\beta<\delta\rangle)
	}. 
\]
Since 
$q\forcestext[\mathbb{P}]{\dot{\xi}_\delta\in \check{F}_\delta}$ 
and by the definition of $H_\delta$, 
we have 
\[
q\forcestext[\mathbb{P}]{
	\check{H}_\delta\subseteq\check{T}_{\dot{\xi}_\delta}
	}.
\]

Now we see 
$q\forcestext[\mathbb{P}]{
	\forall\delta<\gamma\,
	(\check{H}_\delta\subseteq\check{T}_{\dot{\xi}_\delta})
	}$, 
and since $\{H_\delta\st\delta<\gamma\}$ covers $X$, 
we have 
\[
	q\forcestext[\mathbb{P}]{
		\forall\delta<\gamma\,
	(\check{T}_{\dot{\xi}_{\delta}}\in 
	\dot{\sigma}
	(\langle\check{T}_{\dot{\xi}_\beta}\st\beta<\delta\rangle)
	)
	\text{ and }
	\bigcup\{\check{T}_{\dot{\xi}_{\delta}}\st\delta<\gamma\}
	=\check{X}
	}.
\]
This concludes the proof. 
\end{proof}

\begin{rem}
The reader might complain that, 
for $\delta\geq\gamma$, 
$q$ may not force 
	$\check{T}_{\dot{\xi}_{\delta}}$ to be a possible move for \TWO. 
But it is unimportant, 
since the moves after the inning $\gamma$ do not affect the payoff 
and so \TWO may disregard $\check{T}_{\dot{\xi}_{\delta}}$'s 
and take any moves to follow the rule. 
\end{rem}

A similar argument to the above proof 
yields the following corollary. 
An adaptation of the proof 
for the corollary 
is left 
to the reader.

\begin{cor}\label{cor:twowins}
Let $(X,\tau)$ be a 
topological space, 
$\mathbb{P}$ a forcing notion, 
$\alpha$ an ordinal 
and $\kappa=p(X)$. 
If 
\begin{enumerate}
\item \TWO 
	has a winning strategy in 
	$\operatorname{\mathsf{G}}_1^{<\alpha}(\mathcal{O},\mathcal{O})$ 
	on $(X,\tau)$, and  
\item 
	\TWO has a winning strategy in 
	$\operatorname{\mathsf{CG}}^{<\alpha}({<}\kappa)$ on $\mathbb{P}$, 
\end{enumerate}
then 
\[
	\forcestext[\mathbb{P}]{%
		\text{\TWO 
		has a 
		winning strategy in }
		\operatorname{\mathsf{G}}_1^{<\alpha}(\mathcal{O},\mathcal{O})
		\text{ on }
		(\check{X},\tau^{\mathbb{P}})
	}.
\]
\end{cor}

\section{Consequences}\label{sec:conseq}

Theorem~\ref{thm:main} 
together with 
Theorem~\ref{thm:indlindelofgame}
yields the following consequence.

\begin{cor}\label{cor:preserveindlindelof}
Suppose that  
$(X,\tau)$ is an indestructibly Lindel\"{o}f space 
and  
$\mathbb{P}$ is a forcing notion 
such that 
	\TWO has a winning strategy in 
	the game 
	$\operatorname{\mathsf{CG}}^{<\omega_1}({<}\aleph_0)$ on $\mathbb{P}$. 
Then 
\[
	\forcestext[\mathbb{P}]{%
		(\check{X},\tau^{\mathbb{P}})
		\text{ is 
		an indestructibly 
		Lindel\"{o}f space}
	}.
\]
\end{cor}

Now 
consider the following three conditions 
on a Lindel\"{o}f space $(X,\tau)$. 

\begin{enumerate}
\item[(1)]
	$(X,\tau)$  is not 
	destroyed
	by 
	$\operatorname{Fn}(\omega_1,2,\omega_1)$. 
\item[(2)]
	$(X,\tau)$  is indestructibly Lindel\"{o}f, 
	that is, 
	$(X,\tau)$  
	is not 
	destroyed 
	by any ${<}\omega_1$-closed forcing notion.
\item[(3)]
	$(X,\tau)$  is not 
	destroyed
	by any forcing notion $\mathbb{P}$ 
	on which 
	\TWO has a winning strategy in 
	$\operatorname{\mathsf{CG}}^{<\omega_1}({<}\aleph_0)$. 
\end{enumerate}

Clearly $(3)\Rightarrow (2)\Rightarrow (1)$ holds. 
Tall pointed out (see \cite[Theorem~3]{Tall:lindelofpointsgd}) 
a result due to Shelah, 
which claims that 
$(1)\Rightarrow (2)$ holds. 
Corollary~\ref{cor:preserveindlindelof} 
tells us 
that $(2)\Rightarrow (3)$ holds, 
and hence these three conditions are all equivalent. 
In fact, 
we can also give a direct proof of $(1)\Rightarrow (3)$ 
by putting the argument of the proof of Theorem~\ref{thm:main} 
into Shelah's proof.

Using Theorem~\ref{thm:main} 
and 
Theorem~\ref{thm:rothbergergame}, 
we see the following. 

\begin{cor}\label{cor:preserverothberger}
Suppose that  
$(X,\tau)$ has the Rothberger property 
and  
$\mathbb{P}$ is a forcing notion 
such that 
	\TWO has a winning strategy in 
	the game 
	$\operatorname{\mathsf{CG}}^{\omega}({<}\aleph_0)$ on $\mathbb{P}$. 
Then 
\[
	\forcestext[\mathbb{P}]{%
		(\check{X},\tau^{\mathbb{P}})
		\text{ has the Rothberger property}
	}.
\]
\end{cor}

Let $\mathbb{B}(\kappa)$ 
denote the measure algebra 
on $2^{\kappa}$. 
Scheepers and Tall
proved 
that, 
for any infinite cardinal $\kappa$, 
if 
$(X,\tau)$ has the Rothberger property, 
then 
$\forcestext[\mathbb{B}(\kappa)]{
	(\check{X},\tau^{\mathbb{B}(\kappa)})
	\text{ has the Rothberger property}}$ 
\cite[Theorem~15]{ST:lindelof}.  
It is known that 
\TWO has a winning strategy 
in 
the game 
	$\operatorname{\mathsf{CG}}^{\omega}({<}\aleph_0)$ 
	on $\mathbb{B}(\kappa)$ 
\cite{Je:moregame}, 
and hence 
Corollary~\ref{cor:preserverothberger} 
gives an alternate proof of 
their result. 

We also remark that 
Corollary~\ref{cor:preserverothberger} 
extends another result due to 
Scheepers and Tall, 
which claims that 
for a ${<}\omega_1$-closed forcing notion $\mathbb{P}$ 
if 
$(X,\tau)$ has the Rothberger property 
then 
$\forcestext[\mathbb{P}]{
	(\check{X},\tau^{\mathbb{P}})
	\text{ has the Rothberger property}}$ 
\cite[Theorem~21]{ST:lindelof}.

As we mentioned in Section~\ref{sec:pre}, 
a Lindel\"{o}f P-space 
has the Rothberger property. 
Using this fact 
with 
Theorem~\ref{thm:countablechoiceproper}, 
Theorem~\ref{thm:rothbergergame}, 
Proposition~\ref{prop:preservepspace} 
and Theorem~\ref{thm:main}, 
we can deduce the following result. 

\begin{cor}\label{cor:preservelindelofpspace}
Suppose that  
$(X,\tau)$ is a Lindel\"{o}f P-space 
and  
$\mathbb{P}$ is a forcing notion 
such that 
	\TWO has a winning strategy in 
	the game 
	$\operatorname{\mathsf{CG}}^{\omega}(\aleph_0)$ on $\mathbb{P}$. 
Then 
\[
	\forcestext[\mathbb{P}]{%
		(\check{X},\tau^{\mathbb{P}})
		\text{ is a Lindel\"{o}f P-space}
	}.
\]
\end{cor}

Now we can summarize these consequences of 
Theorem~\ref{thm:main} 
as in Table~\ref{table:summary}. 
The table is read as follows: 
``A property of a topological space 
shown in a left-hand column 
is preserved under forcing extension 
by a forcing notion with the property 
shown in the corresponding right-hand column.''

\begin{table}
\caption{Summary of consequences of the main result}
\label{table:summary}
\begin{center}
\begin{tabular}{lll}
\hline
	Topological spaces	
&	Forcing notions
&
\\
\hline
\hline
	(1)~Lindel\"{o}f
&
&
\\
\hline
&	(1)~$\operatorname{Fn}(\omega_1,2,\omega_1)$
&
\\
&	(2)~${<}\omega_1$-closed
&
\\
&	(3)~\TWO has a w.s. in 
		$\operatorname{\mathsf{CG}}^{<\omega_1}(1)$
&
\\
	(2)~indestructibly Lindel\"{o}f
&	(4)~\TWO has a w.s. in 
		$\operatorname{\mathsf{CG}}^{<\omega_1}({<}\aleph_0)$
&	(\ref{cor:preserveindlindelof})
\\
\hline
	(3)~Rothberger
&	(5)~\TWO has a w.s. in 
		$\operatorname{\mathsf{CG}}^{\omega}({<}\aleph_0)$
&	(\ref{cor:preserverothberger})
\\
\hline
	(4)~Lindel\"{o}f P-space
&	(6)~\TWO has a w.s. in 
		$\operatorname{\mathsf{CG}}^{\omega}(\aleph_0)$
&	(\ref{cor:preservelindelofpspace})
\\
\hline
&	(7)~proper
&
\\
\hline
\end{tabular}
\end{center}
\end{table}

Before closing this section, 
we state a consequence 
of Corollary~\ref{cor:twowins}, 
which gives a sufficient condition for a forcing notion 
to preserve the topological property 
``\TWO has a winning strategy in 
$\operatorname{\mathsf{G}}_1^{<\omega_1}(\mathcal{O},\mathcal{O})$ 
on $X$''.

\begin{cor}\label{cor:preservetwowins}
Suppose that 
\TWO has a winning strategy in 
$\operatorname{\mathsf{G}}_1^{<\omega_1}(\mathcal{O},\mathcal{O})$ 
on a 
topological 
space 
$(X,\tau)$ 
and  
$\mathbb{P}$ is a forcing notion 
such that 
	\TWO has a winning strategy in 
	the game 
	$\operatorname{\mathsf{CG}}^{<\omega_1}({<}\aleph_0)$ on $\mathbb{P}$. 
Then 
\[
	\forcestext[\mathbb{P}]{%
		\text{\TWO has a winning strategy in 
		$\operatorname{\mathsf{G}}_1^{<\omega_1}(\mathcal{O},\mathcal{O})$ 
		on }
		(\check{X},\tau^{\mathbb{P}})
	}.
\]
\end{cor}

\section{Discussion}

We will show that, 
under ZFC, 
the assumption 
``\TWO has a winning strategy in 
$\operatorname{\mathsf{CG}}^{\omega}({<}\aleph_0)$ on $\mathbb{P}$'' 
in Corollary~\ref{cor:preserverothberger} 
cannot be weakened to 
``\TWO has a winning strategy in 
$\operatorname{\mathsf{CG}}^{\omega}(\aleph_0)$ on $\mathbb{P}$''.

We use the following famous result due to Laver~\cite{L:borel} 
(also found in \cite[Theorem~8.3.2]{BaJ:set}). 
Let $\mathbb{M}$ denote 
the Mathias forcing notion. 

\begin{thm}\label{thm:mathiassmz}
Suppose that $X$ is an uncountable set of real numbers. 
Then 
\[
\forcestext[\mathbb{M}]{\check{X}\text{ does not have strong measure zero}}. 
\]
\end{thm}

It is easily checked that, 
if 
a forcing notion 
$\mathbb{P}$ 
satisfies Axiom A, 
then 
\TWO has a winning strategy in 
$\operatorname{\mathsf{CG}}^{\omega}(\aleph_0)$ on $\mathbb{P}$. 
On the other hand, 
if
\TWO has a winning strategy in 
$\operatorname{\mathsf{CG}}^{\omega}({<}\aleph_0)$ on $\mathbb{P}$, 
then $\mathbb{P}$ is $\omega^\omega$-bounding 
by Theorem~\ref{thm:finitechoicebounding}. 
The Mathias forcing $\mathbb{M}$ satisfies Axiom A 
but is not $\omega^\omega$-bounding, 
and so  
\TWO has a winning strategy in 
$\operatorname{\mathsf{CG}}^{\omega}(\aleph_0)$ 
on $\mathbb{M}$ 
but 
none
in 
$\operatorname{\mathsf{CG}}^{\omega}({<}\aleph_0)$ 
on $\mathbb{M}$. 

Now assume CH 
and let $L$ be an uncountable Lusin set of real numbers. 
It is known that $L$ has the Rothberger property 
\cite{Ro:propertyc}.
However, 
by Theorem~\ref{thm:mathiassmz} 
we have 
\[
	\forcestext[\mathbb{M}]{
		\check{L}\text{ does not have strong measure zero}}, 
\]
and 
a set of real numbers with the Rothberger property 
has strong measure zero, 
which implies 
\[
	\forcestext[\mathbb{M}]{
		\check{L}\text{ does not have the Rothberger property}}. 
\]

%
%

We do not know 
if the assumption 
``\TWO has a winning strategy in 
$\operatorname{\mathsf{CG}}^{<\omega_1}({<}\aleph_0)$ on $\mathbb{P}$'' 
in Corollary~\ref{cor:preserveindlindelof} 
can be weakened to 
``\TWO has a winning strategy in 
$\operatorname{\mathsf{CG}}^{\omega}({<}\aleph_0)$ on $\mathbb{P}$''. 

\begin{q}
Can we find 
an indestructibly Lindel\"{o}f space $(X,\tau)$ 
and a forcing notion $\mathbb{P}$ 
which satisfy the following?
\begin{enumerate}
\item \TWO has a winning strategy in 
	the game 
	$\operatorname{\mathsf{CG}}^{\omega}({<}\aleph_0)$ on $\mathbb{P}$. 
\item 
	$\forcestext[\mathbb{P}]{
		(\check{X},\tau^{\mathbb{P}})
		\text{ is not a Lindel\"{o}f space}	
		}$.
\end{enumerate}
\end{q}

\begin{rem}
As we mentioned in Section~\ref{sec:conseq}, 
\TWO has a winning strategy in 
$\operatorname{\mathsf{CG}}^{\omega}({<}\aleph_0)$ on 
the measure algebra $\mathbb{B}(\kappa)$ 
for any $\kappa$ 
(moreover, 
for any fixed  $\alpha<\omega_1$, 
\TWO has a winning strategy in 
$\operatorname{\mathsf{CG}}^{<\alpha}({<}\aleph_0)$ 
on $\mathbb{B}(\kappa)$). 
On the other hand, 
it is easy to find a winning strategy for \ONE 
in 
$\operatorname{\mathsf{CG}}^{<\omega_1}({<}\aleph_0)$ 
on $\mathbb{B}(\kappa)$ 
(just note that 
any strictly decreasing sequence of real numbers 
has at most countable order type). 
Unfortunately, 
for any Lindel\"{o}f space $(X,\tau)$ 
we have 
$\forcestext[\mathbb{B}(\kappa)]{
	(\check{X},\tau^{\mathbb{B}(\kappa)})\text{ is Lindel\"{o}f}}$ 
(see \cite{GJT:forcingnormality}). 
\end{rem}

We do not know 
if the assumption 
``\TWO has a winning strategy in 
$\operatorname{\mathsf{CG}}^{\omega}(\aleph_0)$ on $\mathbb{P}$'' 
in Corollary~\ref{cor:preservelindelofpspace} 
can be weakened to 
``$\mathbb{P}$ is proper''.  

Let $\mathbb{CF}$ denote 
the 
poset 
which adjoins a closed unbounded subset of $\omega_1$ 
with finite conditions, 
which is due to Baumgartner \cite{Bau:applpfa}. 
It is known that 
$\mathbb{CF}$
is proper but 
\ONE has a winning strategy in 
$\operatorname{\mathsf{CG}}^{\omega}(\aleph_0)$ on $\mathbb{CF}$. 
So it is natural to ask the following question. 
Note that, 
by Proposition~\ref{prop:preservepspace}, 
a P-space is still a P-space in a forcing extension 
by a proper 
poset. 

\begin{q}\label{q:destroylindelofp}
Is there a Lindel\"{o}f P-space $(X,\tau)$ 
such that 
\[
\forcestext[\mathbb{CF}]{
	(\check{X},\tau^{\mathbb{CF}})\text{ is not Lindel\"{o}f}}\;\text{?}
\]
\end{q}

\subsection*{Acknowledgements}

The author would like to thank Marion Scheepers 
for his helpful comments and discussion during this work.

\noindent
Graduate School of Science\\
Osaka Prefecture University\\
1--1 Gakuen-cho, Naka-ku, Sakai, Osaka 599--8531 JAPAN\\
E-mail: kada@mi.s.osakafu-u.ac.jp

\end{document}